\documentclass{amsart}
\usepackage{amsmath}
\usepackage{amscd}
\usepackage{amssymb}
\usepackage{a4}

\newcommand{\cC}{{\mathcal C}}
\newcommand{\cD}{{\mathcal D}}

\newcommand{\eps}{{\varepsilon}}

\newcommand{\cR}{{\mathcal R}}
\newcommand{\cF}{\mathcal F}
\newcommand{\MA}{{M_H{(A)}}}

\newcommand{\Qml}[1]{Q^l_{max}(#1)}

\newcommand{\Ql}[1]{Q_{cl}(#1)}
\newcommand{\EA}{{\widehat{A}}}
\newcommand{\tA}{{\tilde{A}}}

\newcommand{\End}[1]{{\mathrm{End}}\:(#1)}

\newcommand{\Hom}[2]{{\mathrm{Hom}}\:(#1,#2)}
\newcommand{\EndX}[2]{{{\mathrm{End}}_{#1}\:(#2)}}
\newcommand{\HomX}[3]{{\mathrm{Hom}}_{#1}\:(#2,#3)}

\newcommand{\Ker}[1]{{\mathrm{Ker}}\:(#1)}

\newcommand{\Ann}[2]{{{\mathrm{Ann}_{#1}}\:({#2})}}

\newcommand{\lra}{\longrightarrow}
\newcommand{\ra}{\rightarrow}

\renewcommand{\Im}[1]{{\mathrm{Im}}\:(#1)}
\renewcommand{\lim}{{\displaystyle\varinjlim}}

\newcommand{\AH}{{A\# H}}

\newtheorem*{thm}{Theorem}
\newtheorem*{prop}{Proposition}
\newtheorem*{lem}{Lemma}
\newtheorem*{cor}{Corollary}
\newtheorem*{defn}{Definition}

\begin{document}
\title{When is a smash product semiprime ?}
\date{\today}
%\author{Christian Lomp\\Centro de Matem\'atica, Faculdade de Ci\^{e}ncias,\\
%Universidade do Porto, R.Campo Alegre 687,\\
%4169-007 Porto, Portugal\\
%clomp@fc.up.pt}
\author{Christian Lomp}
\address{Centro de Matem\'atica, Faculdade de Ci\^{e}ncias,\\ Universidade do Porto, R.Campo Alegre 687,\\ 4169-007 Porto, Portugal}
%\email{clomp@fc.up.pt}

\thanks{Work supported by {\it Funda{\c{c}\~ao} para a Ci\^encia e a Tecnologia}
(FCT) through the {\it Centro de Matem\'atica da Universidade do
Porto } (CMUP)}

\begin{abstract}
It is an open question whether the smash product of a semisimple
Hopf algebra and a semiprime module algebra is semiprime. In this
paper we show that the smash product of a commutative semiprime
module algebra over a semisimple cosemisimple Hopf algebra is
semiprime. In particular we show that the central $H$-invariant
elements of the Martindale ring of quotients of a module algebra
form a von Neumann regular and self-injective ring whenever $A$ is
semiprime. For a semiprime Goldie PI $H$-module algebra $A$ with
central invariants we show that $\AH$ is semiprime if and only if
the $H$-action can be extended to the classical ring of quotients
of $A$ if and only if every non-trivial $H$-stable ideal of $A$
contains a non-zero $H$-invariant element. In the last section we
show that the class of strongly semisimple Hopf algebras is closed
under taking Drinfeld twists.  Applying some recent results of
Etingof and Gelaki we conclude that every semisimple cosemisimple
triangular Hopf algebra over a field is strongly semisimple.
\end{abstract}
 \maketitle
 %\tableofcontents

\section{Introduction}
It is an important open question in the theory of Hopf algebra
actions whether the smash product $\AH$ of a semisimple Hopf $H$
and a semiprime left $H$-module algebra $A$, is semiprime (see
\cite[Question 4.4.7]{Montgomery}).

Fisher and Montgomery had proved an analogous result for group
rings (see \cite{FisherMontgomery}) and Cohen and Montgomery for
duals of group rings (see \cite{CohenMontgomery}). Attempts had
been made to tackle this question often by restricting the class
of Hopf algebras. ( see for example \cite{MontgomerySchneider99}).

In order to give a partial answer to the semiprimness question, we
will restrict the class of module algebras rather than the class
of Hopf algebras. In particular we will show that the question has
a positive answer for commutative module algebras in
characteristic $0$. The main step is to show that the subring of
central $H$-invariant elements of the Martindale ring of quotients
is von Neumann regular. The result follows applying a theorem of
S.Zhu which says that a commutative module algebra is an integral
extension of its invariants if the Hopf algebra involved is
semisimple and cosemisimple.

In general one might ask what are necessary or sufficient
conditions for a smash product to be semiprime. A very important
necessary condition is the existence of non-trivial $H$-invariant
elements in non-zero $H$-stable ideals of the module algebra. A
sufficient condition is the ability of extending the $H$-action on
a semiprime Goldie module algebra to its classical ring of
quotients. We will see in Theorem \ref{SemiprimeGoldieZentral}
that for semiprime Goldie PI module algebras with central
invariants those conditions are equivalent to the smash product
being semiprime. In the final section we show that the class of
strongly semisimple Hopf algebras is closed under Drinfeld twists.
Applying finally a recent result of Etingof and Gelaki, we can
also conclude that triangular semisimple cosemisimple Hopf
algebras are strongly semisimple and satisfy the property that
their smash product with a semiprime module algebra is semiprime.

All rings are supposed to be associative and have a unit element
unless otherwise stated. Throughout the text
$R$ will denote a commutative ring, $H$ a Hopf algebra over $R$ with antipode $S$, counit $\eps$ and comultiplication $\Delta$.
We will make use of the so-called Sweedler-notation $\Delta(h)=\sum_{(h)} h_1 \otimes h_2$ for the comultiplication
of an $h\in H$. A left $H$-module algebra $A$ is an $R$-algebra in the category of left $H$-modules.
The smash product of $A$ and $H$ is an $R$-algebra with underlying $R$-module $A{\otimes_R} H$ and denoted by $\AH$.
The multiplication of two elements $a\# h$ and $b\# g$ in $\AH$ is defined to be equal to
$$ (a\# h)(b\# g) := \sum_{(h)} a(h_1 b) \# h_2g.$$
We emphasis that $A$ is a cyclic left $\AH$-module and $\EndX{\AH}{A} \simeq A^H$. This allows to study
$A$, $A^H$ and $\AH$ in module-theoretic terms.

We refer to all unexplained Hopf-algebraic terms to
\cite{Montgomery} and \cite{Pareigis73V}, to all ring-theoretic
terms to \cite{Lam} and to all module-theoretic terms to
\cite{wisbauer}.

\section{Separability of smash products}

Many results on group actions are stated in terms of algebras over
rings rather than in terms of algebras over fields.
Throughout the paper we will consider Hopf algebras over a commutative ring $R$.
Just when applying deeper results on Hopf algebras over fields we
will assume that $R$ is a field. In the case of a base ring $R$
the adequate analogue of a semisimple Hopf algebra (over a field) is
a Hopf algebra that is separable over $R$.
\subsection{}We will shortly recall
the definition of separability in non-commutative ring extensions (see \cite{HirataSugano}).
\begin{defn} Let $S\subseteq T$ be any ring
extension. $T$ is called {\it separable} over $S$ if there exists
an idempotent
$$\omega:= \sum_{i=1}^n x_i\otimes y_i \in T {\otimes_S} T \mbox{ such that } \sum_{i=1}^nx_iy_i = 1 \mbox{ and  } t\omega = \omega t$$
holds for all $t\in T$.
We refer to $\omega$ as the {\it separability idempotent} of $T$ over $S$.
\end{defn}

Here we consider $T{\otimes_S} T$ as a $T-T$-bimodule via
$t(x\otimes y) = tx \otimes y$ and $(x\otimes y)t = x \otimes yt$ for all $t\in T$ and $x\otimes y \in T {\otimes_S} T$.

\subsection{}
Separable extensions are in particular semisimple extensions (see
\cite{HirataSugano}). An extension $S\subseteq T$ is called semisimple
if every exact sequence of left $T$-modules, which splits as a
sequence of left $S$-modules, splits. Hence if $H$ is a Hopf
algebra over some field $k$ such that $k\subseteq H$ is separable,
$H$ must be a semisimple ring. (Note that `semisimple ring' shall always mean
`semisimple artinian ring'). We will see soon that the converse is true
as well.

\subsection{}\label{Separability}
Recall the submodule of left integrals in a Hopf algebra $H$:
$$\int_l :=\{ t\in H \mid \forall h\in H: ht=\eps(h)t \}.$$
Right integrals are defined analogously.
It is known that $\int_l\neq 0$ in case $H$ is finitely generated
and projective as $R$-module (see \cite{Pareigis71}).
The following Proposition gives a criterium for $\AH$ to be
separable over $A$.
\begin{prop}
Let $H$ be a Hopf algebra over $R$ and let $A$ be a left
$H$-module algebra. Assume that there exists a left or right
integral $t$ in $H$ with $\varepsilon(t)1_A$ invertible in $A$.
Then $\AH$ is separable over $A$.
\end{prop}
\begin{proof}
Let $t$ be a right integral in $H$ such that $\eps(t)1_A$ is
invertible in $A$ and let $z\in A$ be its inverse. For any $a\in
A$:
$$ (az-za) = (az-za)\eps(t)z = (a-a)z=0$$
implies $z\in Z(A)$ and for any $h\in H$:
$$ (h\cdot z - \eps(h)z) = (h\cdot z -
\eps(h)z)\eps(t)z = (h\cdot 1_A - \eps(h)1_A)z=0$$ shows $z\in
A^H$. Hence $z\in Z(A)^H:=Z(A)\cap A^H$. Consider the element
$$\omega := \sum_{(t)} [1\# S(t_1)] \otimes [z\# t_2] \in \AH \otimes_A \AH.$$
We will show that $\omega$ is a separability idempotent for $\AH$.

Let $\mu:\AH\otimes_A \AH\longrightarrow \AH$ denote the
multiplication map. We have
\begin{eqnarray*}
\mu(\omega) = \sum_{(t)} (1\# S(t_1))(z\# t_2) &=& \sum_{(t)} S(t_2)\cdot z \# S(t_1)t_3\\
&=& \sum_{(t)} z \# S(t_1)t_2 = \eps(t)z \# 1_H = 1_A \# 1_H.
\end{eqnarray*} Let
$a\in A$. Then the following holds:
\begin{eqnarray*}
\omega a&=& \sum_{(t)} (1\# S(t_1))\otimes (z\# t_2)(a\#1)\\
&=& \sum_{(t)} (1\# S(t_1)) \otimes (z (t_2 a) ) \# t_3)\\
&=& \sum_{(t)} (1\# S(t_1))(t_2 a \# 1)\otimes (z\# t_3) \\
&=& \sum_{(t)} (S(t_2)t_3 a  \# S(t_1))\otimes (z\# t_4) \\
&=& \sum_{(t)} (a \# S(t_1) )\otimes (z\# t_2) = a \omega
\end{eqnarray*}
Hence $\omega a = a \omega$ shows that $\omega$ is
$A$-centralising. Let $h\in H$ and note that
$$ h\otimes \Delta(t) = \sum_{(h)} h_1 \otimes \Delta(t\eps(h_2))
= \sum_{(h)} h_1 \otimes \Delta(th_2) = \sum_{(h,t)} h_1 \otimes
t_1h_2 \otimes t_2h_3$$
holds. Using this equation we get:
\begin{eqnarray*}
h \omega &=& \sum_{(t)} (1\# hS(t_1))\otimes (z\# t_2)\\
&=& \sum_{(h,t)} (1\# h_1S(t_1h_2)) \otimes (z \# t_2h_3)\\
&=& \sum_{(t)} (1\# h_1S(h_2)S(t_1))\otimes (z\# t_2 h_3) \\
&=& \sum_{(t)} (1\# S(t_1))\otimes (z\# t_2h) = \omega h
\end{eqnarray*}
showing that $\omega$ is $H$-centralising. Thus $\omega$ is a
separable idempotent of $\AH$ over $A$. For a left integral $t$
with $\eps(t)$ invertible in $A$ we set $t':=S(t)$. Since $t'$ is
a right integral and $\eps(t')=\eps(t)$ we can argue as above and
conclude that $\AH$ is separable over $A$.
\end{proof}

\subsection{}
Letting $H$ act trivially on $R$ by setting $hr:=\eps(h)r$ for all $h\in H$ and $r\in R$, $R$ becomes
a left $H$-module algebra and $R\# H \simeq H$. Proposition \ref{Separability} shows that
$H$ is separable over $R$ if and only if there exists a left or right integral $t$ in $H$ with $\eps(t)$ invertible in $R$.
While the sufficiency follows from the Proposition, the necessity follows because if $H$ is separable over $R$ the
$H$-linear map $\eps: H \rightarrow R$ splits as right $H$-modules. Hence there exists an $H$-linear
$\sigma:R \rightarrow H$ such that $\eps(\sigma(1))=1$. The element $t:=\sigma(1)$ is our right integral.

In particular $\AH$ is separable over $A$ for every left $H$-module algebra $A$ whenever $H$ is separable over $R$
Thus for instance whenever $H$ is a semisimple Hopf algebra over a field.

Note that this fact holds without assuming any additional hypothesis
on $H$ as a module over $R$. On the other hand it is well-known that a separable
$R$-algebra $H$ must be finitely generated as $R$-module if $H$
is projective as $R$-module.

\subsection{}
In case of a group ring $H=R[G]$ with $G$ a finite group. The submodule of left and right
integrals $\int_l$ is spanned by the element $t:=\sum_{g\in G} g$. For an $R$-algebra $A$ where $G$ acts on, $A\# G$ is
equal to the skew group ring of $A$ and $G$. Proposition \ref{Separability} says that $A\# G$ is separable over $A$ provided
$\eps(t)=|G|$ is invertible in $A$.

\subsection{}Since separable extensions $S\subseteq T$ are semisimple extensions, every left $T$-module that is
projective as left $S$-module is also projective as left $T$-module.
In particular any separable extension of a semisimple artinian ring is itself semisimple artinian.
Our next Lemma shows that the analogue statement for flat modules and von Neumann regular rings
is also true.

\begin{lem}\label{LemmaFlatSeparable} Suppose $T$ is separable over a subring $S$. Then every left
$T$-module that is flat as left $S$-module is also flat as left $T$-module.
\end{lem}

\begin{proof}
Let $M$ be a left $T$-module such that $M$ is flat as left
$S$-module and let $\beta: T\otimes_S M \longrightarrow M$ with
$\beta(t\otimes m):=tm$. Obviously $\beta$ is an epimorphism of
left $T$-modules. Consider the embedding $i: S \longrightarrow T$
in Mod-$S$. Since $M$ is flat as left $S$-module we get an
embedding: $\alpha:=i\otimes 1_M : M=S\otimes_S M \longrightarrow
T\otimes_S M$ where $\alpha(m):=1_T \otimes m$. Hence $\alpha$
lets $\beta$ split as $S$-module homomorphisms. Since $T$ is
separable over $S$  there exists also a left $T$-module map
$\alpha'$ that lets $\beta$ split, i.e. $M$ is a direct summand of
$T\otimes_S M$ as left $T$-module. Since $M$ is flat as left
$S$-module, $M$ is isomorphic to the direct limit of some finitely
generated projective left $S$-modules $P_\lambda$, i.e. $ M \simeq
\lim P_\lambda.$ Hence $T\otimes_S M \simeq \lim T\otimes_S
P_\lambda$ is also a direct limit of finitely generated projective
left $T$-modules $T\otimes_S P_\lambda$, i.e. $T\otimes_S M$ is a
flat left $T$-module. As a direct
summand of a flat $T$-module, $M$ is also flat as $T$-module.
\end{proof}

Hence a separable extension of a von Neumann regular ring is itself von Neumann regular.

\subsection{} Combining Lemma \ref{LemmaFlatSeparable} and Proposition \ref{Separability} we get the following
important Corollary which generalises a result of Cohen and Fischman that says that
$\AH$ is semisimple whenever $A$ and $H$ are semisimple (see \cite[Theorem 6]{CohenFishman}).
\begin{cor}\label{VonNeumannRegularesSmashProdukt}
Let $H$ be an $R$-Hopf algebra and $A$ a left $H$-module algebra,
such that there exists a left or right integral $t$ in $H$ with
$\eps(t)$ invertible in $A$. If $A$ is von Neumann regular, then $\AH$ is von  Neumann regular.
If $A$ is semisimple artinian, then $\AH$ is semisimple artinian.
\end{cor}
\subsection{}
A first application of the corollary above will allow us to show
that whenever the $H$-action can be extended to the left maximal
ring of quotients $\Qml{A}$ of a left non-singular $H$-module
algebra $A$ the smash product $\AH$ must also be left non-singular
and moreover its left maximal ring of quotients is isomorphic to
$\Qml{A} \# H$.

Recall the definition of the left maximal ring of quotients. Let
$S$ be any ring and denote by $E(S)$ its injective hull in
$S$-Mod. Define the left maximal ring of quotients of $S$ as the
$S$-submodule $$\Qml{S} := \{ m \in E(S) \mid \forall f \in \EndX{S}{E(S)}: f(S)=0 \Rightarrow f(m)=0 \}$$
of $E(S)$.
Let $B:=\EndX{\EndX{S}{E(S)}}{E(S)}$ be the biendomorphism ring of
$E(S)$. The evaluation map $\Psi: B \rightarrow \Qml{S}$ with
$\varphi \mapsto \varphi(1)$ is an isomorphism of abelian groups
and induces a ring structure on $\Qml{S}$. Hence one might
identify $\Qml{S}$ with the biendomorphism ring of the injective
hull of $S$.

\subsection{} Recall that a submodule $N$ of a module $M$
is called {\it dense} whenever $\HomX{}{L/N}{M}=0$ for all $N\subseteq L \subseteq M$. $\Qml{S}$
can also be seen as the maximal extension $E$ of $S$ such that $S$
is dense in $E$.

\begin{lem}\label{LemmaForNichtSingular}
Let $S\subseteq T$ be a ring extension such that $\HomX{S}{T/S}{T}=0$ and $_ST$ is
injective. Then $T\simeq \Qml{S}$ as rings.
\end{lem}

\begin{proof}
Let $L$ be an $S$-submodule of $T$ containing $S$. By injectivity
of $T$, every homomorphism $f:L/S\rightarrow T$ can be extended to
an homomorphism $\bar{f}:T/S\rightarrow T$ which is zero by
hypothesis. Thus $S$ is dense in $T$. By \cite[13.11]{Lam} there
exists an injective ring homomorphism $g:T \hookrightarrow
\Qml{S}$ such that $g(s)=s$ for all $s\in S$. Hence $g$ is left
$S$-linear and by injectivity of $T$, $\Im{g}$ is a direct summand
of $\Qml{S}$ containing the essential submodule $S$. Thus $g$ must
be surjective and must be an isomorphism of rings.
\end{proof}

\subsection{} In the following Theorem  we will apply Corollary
\ref{VonNeumannRegularesSmashProdukt} and Lemma
\ref{LemmaForNichtSingular} to show that $\Qml{\AH}\simeq
\Qml{A}\# H$ is von Neumann regular. Using Johnson's Theorem that
states that a ring $S$ is left non-singular if and only if its
left maximal ring of quotients $\Qml{S}$ is von Neumann regular we
can conclude that $\AH$ is left non-singular.

\begin{thm}\label{nichtsingulaer} Let $H$ be a Hopf algebra
over $R$ with $H_R$ finitely generated and projective.
Let $A$ be a left $H$-module algebra, such that there exists a
left or right integral $t\in H$ with $\eps(t)1_A$ invertible in $A$.
Assume that the $H$-action extends to the left maximal ring of
quotients $\Qml{A}$. If $A$ is left non-singular, then $\AH$ is
left non-singular and $\Qml{\AH}\simeq \Qml{A}\# H$.
\end{thm}

\begin{proof} By hypothesis $A$ is left non-singular.
Hence by Johnson's Theorem the maximal ring of quotients
$Q:=\Qml{A}$ of $A$ is von Neumann regular and equals $E(A)$
the injective hull of $A$. In particular $Q$ is injective as $A$-module.
The invertibility of $\eps(t)$ in $A$ (and hence in $Q$) implies
the separability of $Q\# H$ over $Q$ by Proposition \ref{Separability}.
 From Corollary \ref{VonNeumannRegularesSmashProdukt} we know that
 $Q\# H$ is von Neumann regular.
Applying the exact functor $-\otimes_R H$ to the exact sequence
$$\begin{CD} 0@>>> A @>>> Q @>>> Q/A @>>> 0 \end{CD}$$
we get $Q\# H/\AH \simeq (Q/A)\otimes_R H $ as left $A$-modules.
Since $_RH$ is a direct summand of a free module $R^k$ with $k\geq
1$ and since $A$ is dense in $Q$, we get:
\begin{eqnarray*}
\HomX{\AH}{Q\# H/\AH}{Q\# H} &\subseteq&
\HomX{A-}{(Q/A)\otimes_R H}{Q\otimes_R H} \\
&\subseteq& \HomX{A-}{(Q/A)^k}{Q^k} = 0.
\end{eqnarray*}
Hence $\HomX{\AH}{Q\# H/\AH}{Q\# H}=0$. Since $\AH$ is separable
over $A$, we can also conclude that $Q\# H$ is an injective left
$\AH$-module, as $Q$ and $Q\# H$ are injective left $A$-modules.
By Lemma \ref{LemmaForNichtSingular} $\Qml{\AH} \simeq Q \# H$
and by Johnson's Theorem (see \cite[13.36]{Lam}) $\AH$ is left
non-singular.
\end{proof}

\subsection{}
The question whether the $H$-action of a semisimple Hopf algebra can be extended to the maximal ring of
quotients of a module algebra is still open. A claim that this is always possible was made in
\cite{Rumynin} but its proof is not complete as it was confirmed by the author of \cite{Rumynin}.

\section{Commutative semiprime module algebras}

Consider the subring $\MA$ of $\EndX{R}{A}$ generated by the $H$-action
on $A$ and by the left and right multiplications of elements of
$A$:
$$\MA := \left< \{ L_a, R_a, L_h \mid a \in A, h\in H \}\right>
\subseteq \EndX{R}{A},$$ where $L_a$ and $R_a$  denotes the left
and right multiplication with $a\in A$, respectively, and $L_h$
denotes the $H$-action of the element $h$ on $A$. $A$ is a cyclic
faithful $\MA$-module whose submodules are precisely the
$H$-stable two-sided ideals of $A$. If $A$ is commutative then $\MA \simeq
\AH/\Ann{\AH}{A}$.

\subsection{}A module algebra $A$ is called {\it $H$-semiprime} if $A$ does not
contain any non-trivial nilpotent $H$-stable ideals.

\begin{lem}\label{LemmaRational1}
The following statements are equivalent for an $H$-stable ideal
$I$ of an $H$-semiprime module algebra $A$.
\begin{enumerate}
        \item[(a)] $l.ann_A(I)=0$;
        \item[(b)] $I$ is an essential $\MA$-submodule of $A$;
        \item[(c)] $I$ is a dense $\MA$-submodule of $A$.
\end{enumerate}
\end{lem}

\begin{proof}
(a) $\Rightarrow$ (b) Let $J$ be an $H$-stable ideal of $A$. Since
the left annihilator of $I$ is zero, $J\cap I \supseteq JI \neq 0$ shows that $I$ is essential $\MA$-submodule of $A$.\\
(b) $\Rightarrow$ (c) Let $J$ be an $H$-stable ideal of $A$ containing $I$ and let
$f:J\lra A$ be $\MA$-linear such that $I\subseteq \Ker{f}$. Then $K:= f(J)\cap
I$ is nilpotent since $K^2 \subseteq f(J)I = f(I)=0$. As $A$ is
$H$-semiprime  $K=0$ and as $I$ is essential $f=0$. Hence $\HomX{\MA}{J/I}{A}=0$ shows that $I$ is dense in $A$.\\
(c) $\Rightarrow$ (a) Let $J$ denote the left annihilator of $I$.
Since for all $h\in H, x\in J$ and $y\in I$ we have:
$$(hx)y = \sum_{(h)} h_1( x (S(h_2)y)) = 0,$$
$J$ is an $H$-stable ideal of $A$. Since $A$ is $H$-semiprime,
$I\cap J=0$. Let $\pi:J\oplus I \lra J$ be the projection, then
$\pi\in\HomX{\MA}{(J\oplus I)/I}{A}=0$. Hence $I$ has zero left annihilator.
\end{proof}

\subsection{}Recall that the self-injective hull $\widehat{M}$ of a module $M$
is the largest $M$-generated submodule of its injective hull
$E(M)$. The endomorphism of the self-injective hull of a module
whose essential submodules are dense is known to be von Neumann
regular and self injective (see \cite[11.2]{wisbauer96}). Applying
this module-theoretic fact to our situation Lemma \ref{LemmaRational1} shows that the endomorphism
ring $T$ of the self-injective hull $\EA$ of $A$ as $\MA$-module is von Neumann regular and self-injective.
We will construct an isomorphism between $T$ and the subring of central $H$-invariant elements of the Martindale
ring of quotients of $A$.

\subsection{}Let $\cF$ denote the set of ideals of $A$ with zero left and right
annihilator. The {\it right Martindale ring of quotients} of $A$ is
$$ Q(A):=\lim \{ \HomX{-A}{I}{A} \mid I\in \cF\}.$$
Alternatively one might construct $Q(A)$ as follows: define an
equivalence relation on $\bigcup_{I\in \cF} \HomX{-A}{I}{A}$ by
letting $f:I\lra A$ to be equivalent to $g:J\lra A$ if there exists a
$K\in \cF$ such that $K\subseteq I\cap J$ and $f_{\mid K} =
g_{\mid K}$. Note that the equivalence class of the zero map
contains all maps $f$ that vanish on some ideal in $\cF$. Addition
is defined by $[f]+[g]:=[f+g: I\cap J \lra A]$ while multiplication is
set to be $[f][g]:=[fg: JI \lra A]$ where $fg$ denotes the composition map $a \mapsto f(g(a))$.

In order to extend the $H$-action on $A$ to some subring of
$Q(A)$, Miriam Cohen considered the subset $\cF_H$ of $H$-stable
ideals belonging to $\cF$ and constructed the following ring:
$$Q_0(A):= \lim \{\HomX{-A}{I}{A} \mid I\in \cF_H\}.$$
We will refer to the elements of $Q_0(A)$ as equivalence classes
in the above sense. Moreover $Q_0(A)$ is a subring of $Q(A)$. The
$H$-action on $A$ extends to $Q_0(A)$ by letting an
element $h\in H$ act on $f:I\lra A$ by $(h\cdot f):I\lra A$ with
$$(h\cdot f)(x) := \sum_{(h)} h_1 f( S(h_2)x ) \:\mbox{ for all } x\in
I.$$ One checks as in \cite[Theorem 18]{cohen85} that $Q_0(A)$ becomes a left
$H$-module algebra with this action.

\subsection{}
We are now in position to show that the subring of central
$H$-invariant elements $Z(Q_0)^H:=Z(Q_0)\cap Q_0^H$ of the right Martindale ring of
quotients of a semiprime module algebra is von Neumann regular and
self-injective.

\begin{prop}\label{ZenterOfMartindale} Let $H$ be a Hopf algebra over $R$ and let $A$ be a left
$H$-semiprime module algebra with right Martindale ring of
quotients $Q_0$. Let $T$ be the endomorphism ring of the
self-injective hull $\EA$ of $A$  as $\MA$-module. Assume that $A$
is commutative or $A$ is semiprime or $H$ has a bijective
antipode. Then
$$\psi: T \rightarrow Z(Q_0)^H \mbox{ with } f \mapsto [f : I_f \rightarrow A]$$
is a ring isomorphism where $I_f := f^{-1}(A)\cap A$.
Moreover $Z(Q_0)^H$ is a von Neumann regular self-injective ring.
\end{prop}

\begin{proof}
Let $\EA$ denote the self-injective hull of $A$ as $\MA$-module
and let $T$ denote the endomorphism ring of $\EA$ as $\MA$-module.
For each endomorphism $f\in T$ define $I_f := f^{-1}(A)\cap A.$ Since pre-images of
essential submodules are essential, $I_f$ is an essential
$\MA$-submodule of $A$. By Lemma \ref{LemmaRational1} $I_f$ has zero left annihilator.
If $A$ is commutative or semiprime
$I_f$ has also zero right annihilator and belongs to $\cF_H$. If the antipode of $H$ is bijective then the
right annihilator $J$ of $I_f$ is also an $H$-stable ideal since
for all $h\in H, x\in J$ and $y\in I_f$ we have:
$$x(hy) = \sum_{(h)} h_2( (S^{-1}(h_1) x) y) = 0.$$
As $I_f \cap J$ is a nilpotent $H$-stable ideal and as $A$ is $H$-semiprime $I_f \cap J$
must be equal to the zero submodule.
$I_f$ being an essential $\MA$-submodule implies that $J$ is zero. Thus also in this case $I_f$ belongs to $\cF_H$.

We will show that $\psi$ is a ring homomorphism. Let $f,g \in T$.
Note that $I_fI_g \in \cF_H$ and $I_fI_g \subseteq I_{fg}$
%since
%for all $xy \in J$ the following holds :
%$$fg(xy)=f(xg(y))=f(x)g(y)\in A.$$
Thus
$$\psi(f)\psi(g) = [f:I_f\ra A][g:I_g \ra A]
= [fg: I_fI_g \ra A] = [fg:I_{fg}\ra A]=\psi(fg).$$ This shows that
$\psi$ is a ring homomorphism. Assume $\psi(f)=0$ for some $f\in
T$. Then there exists an $J\in \cF_H$ with $J\subseteq I_f$ and
$f(J)=0$. Hence $f\in \HomX{\MA}{I_f/J}{A}=0$ as $J$ is dense by
Lemma \ref{LemmaRational1}. This shows that $\psi$ is injective.
On the other hand $\psi$ is also surjective. Let $[q:I \lra A] \in
Z(Q_0)^H$. First note that $q$ is $\MA$-linear: Let $a\in A$ then
$[q][L_a]=[L_a][q]$ implies the existence of an ideal $J\in\cF_H$
with $J\subseteq I$ and
$$q':=qL_a-L_a q \in \HomX{-A}{I/J}{A}.$$
Since $J$ has zero left annihilator and $q'(I/J)J=0$ we can
conclude $q'=0$. This shows
$$q(ax) = qL_a(x) = L_aq(x) = aq(x)$$
for all $x \in I$. Hence $q$ is a left $A$-linear.

Note that since $q\in Q_0^H$ for all $h\in H: h\cdot q = \eps(h)q$.
Let $h\in H$. For all $x\in I$ we have:
\begin{eqnarray*}
q(hx) = \sum_{(h)} \eps(h_1)q(h_2x)
&=& \sum_{(h)} h_1 \cdot q(h_2 x)\\
&=& \sum_{(h)} h_1 q( S(h_2)h_3 x)
= \left[\sum_{(h)} h_1\eps(h_2)\right] q(x) = h q(x).
\end{eqnarray*}
This shows the $H$-linearity of $q$. Since $q$ is by definition right $A$-linear we have shown that
$q$ is an $\MA$-linear map.

By injectivity of $\EA$, $q:I\rightarrow A$ can be extended to an
$\MA$-linear map $\bar{q}\in T$. This extension is unique since
$\HomX{\MA}{\EA/I}{\EA}=0$. Moreover
$\psi(\bar{q})=[q]$ as $I\subseteq I_{\bar{q}}$ and $\bar{q}_{\mid
I}=q$. This shows that $\psi$ is surjective and we have
established an isomorphism of rings between $Z(Q_0)^H$ and $T$
which is von Neumann regular and self-injective by \cite[11.2]{wisbauer96}.
\end{proof}

\subsection{} Our main result follows now easily from the preceding paragraphs.

\begin{thm}\label{CohenLsgCommutative}
Let $H$ be a Hopf algebra over $R$ such that $H_R$ is flat and let $A$ be a commutative semiprime left $H$-module
algebra. Assume that there exists a left or right integral $0\neq t\in H$ such that $\eps(t)$ is not a
zero divisor in $A$.
Then $\AH$ is semiprime provided $A$ is integral over $A^H$.
\end{thm}

\begin{proof}
Denote by $Q_0$ the right Martindale ring of quotients of the module algebra $A$.
Assume $\eps(t)$ is invertible in $A$. Let $\tA:=<A,Q_0^H> \subseteq
Q_0$ be the subalgebra of $Q_0$ generated by $A$ and $Q_0^H$. Obviously
$\tA$ is a left $H$-module algebra. Since $\tA$ is a subalgebra of the right Martindale ring of quotients $Q$ of $A$ which is
commutative and semiprime, also $\tA$ is commutative and
semiprime. By hypothesis $A$ is an integral extension of
$A^H$. Hence $\tA$ is integral over $Q_0^H$. To see this note that
$A^H \subseteq Q_0^H$ and let $aq \in \tA$. There exists a monic
polynomial
$$f(X)=\sum_{i=0}^n r_iX^i \in A^H[X]$$ with $f(a)=0$.
Define the monic polynomial
$$\tilde{f}(X):=\sum_{i=0}^n r_iq^{n-i} X^i \in Q_0^H[X].$$
Then $\tilde{f}(aq) = f(a)q^n = 0$ shows that every element of the form
$aq$ of $\tA$ is integral over $Q_0^H$. Since the set of integral
elements is closed under sums, we get $\tA$ is integral over
$Q_0^H$.

By Proposition \ref{ZenterOfMartindale} $Q_0^H$ is von Neumann regular. Recall
that a commutative ring is von Neumann regular if and only if it
is semiprime and every prime ideal is maximal. Since $Q_0^H\subseteq
\tA$ is an integral extension, the height of a prime ideal $P$ in
$\tA$ is equal to the height of the prime ideal $P\cap Q_0^H$ (see for example \cite[9.2]{Eisenbud})
every prime ideal of $\tA$ is maximal and therefore $\tA$ is von Neumann regular.

Since $\eps(t)1_A$ is invertible in $A$, it is also invertible in
$\tA$. By Corollary \ref{VonNeumannRegularesSmashProdukt} $\tA\# H$ is von
Neumann regular. Let $I\subseteq \AH$ be an ideal with $I^2=0$.
Then $\tilde{I}:=I(Q_0^H\# 1)$ is an ideal of $\tA\# H$. Since
$Q_0^H\# 1$ is central in $\tA\# H$ we get $\tilde{I}^2=0$. As
$\tA\# H$ is von Neumann regular, hence semiprime, we have
$\tilde{I}=0$. Since $_RH$ is flat, $\AH$ is a subring of $\tA \# H$ and thus $I=0$.
This shows $\AH$ does not  contain a non-trivial nilpotent ideal and must be semiprime.\\
In case $\eps(t)1_A$ is not invertible in $A$ but a non-zero divisor, we can localise $A$ by
the powers of $\eps(t)1_A$ and obtain a semiprime commutative
module algebra $A[\eps(t)^{-1}]$. Thus $A[\eps(t)^{-1}] \# H =
\AH[\eps(t)^{-1}\# 1]$ is semiprime and so must be also $\AH$.
\end{proof}

\subsection{} S.Zhu showed that a commutative $H$-module algebra $A$ is an integral extensions of its invariants whenever
$H$ is a finite dimensional Hopf algebra over a field $k$ such
that $char(k)\nmid dim(H)$ and $S^2=id$ (see \cite[Theorem
2.1]{Zhu96}). Etingof and Gelaki proved in \cite{EtingofGelaki}
that a finite dimensional Hopf algebra $H$ satisfies $char(k)\nmid
dim(H)$ and $S^2=id$ if and only if $H$ is semisimple and
cosemisimple. Combining Zhu's and Etingof and Gelaki's result with
Theorem \ref{CohenLsgCommutative} we obtain the following
Corollary.

\begin{cor}\label{CommutativeCase} Let $H$ be a semisimple cosemisimple Hopf algebra over a field
and let $A$ be a commutative semiprime $H$-module algebra. Then $\AH$ is semiprime.
\end{cor}

It is well known, that a semisimple Hopf algebra  over a field of
characteristic $0$ is also cosemisimple.

\section{Semiprime Goldie PI Module algebras}

Assume that the smash product $\AH$ of a module algebra $A$ and a
semisimple Hopf algebra $H$ is semiprime. Then every non-zero
$H$-stable left ideal of $A$ contains a non-zero $H$-invariant
element. In this section we will show that this necessary
condition is also a sufficient condition for semiprime Goldie PI
module algebras with central invariants. More generally we will
show that the $H$-action on such a module algebra can be extended
to its classical ring of quotients in case every non-zero
$H$-stable left ideal contains a non-zero $H$-invariant element.

\subsection{} A module $M$ is called {\it retractable} if $\Hom{M}{N}
\neq 0$ for all non-zero submodules $N$ of $M$ (see
\cite{Zelmanowitz96}). Recall that one has an $R$-linear
isomorphism $I^H \simeq \HomX{\AH}{A}{I}$ for all $H$-stable left
ideals $I$ of $A$. Hence the existence of non-trivial
$H$-invariant elements in non-zero $H$-stable left ideals can be
expressed as $A$ being a retractable $\AH$-module.

\begin{lem}\label{RetractableSemisimple} Let $M$ be a retractable left $R$-module whose
endomorphism ring is semisimple. Then  $M$ is a semisimple
artinian $R$-module. If moreover $R$ is a PI-ring, then $M$ is
finitely generated over its endomorphism ring.
\end{lem}

\begin{proof}
Let $N$ be a non-zero submodule of $M$. By hypothesis there exists
a non-trivial idempotent $e\in S:=\EndX{R}{M}$ such that
$\HomX{R}{M}{N} = Se$. Thus $M=Me \oplus M(1-e)$ implies $N=Ne
\oplus (N \cap M(1-e))$.
%Note that $\HomX{R}{M}{M(1-e}=S(1-e)$.
Hence $$\HomX{R}{M}{N\cap M(1-e)} = \HomX{R}{M}{N}\cap
\HomX{R}{M}{M(1-e)} = Se \cap S(1-e)=0$$ implies by hypothesis
$N\cap M(1-e)=0$, i.e. $N$ is a direct summand of $M$. This shows
that $M$ is a semisimple $R$-module. As $\End{M}$ is artinian, $M$ is artinian.\\
Write $M=\oplus_{i=1}^k E_i^{n_i}$ with pairwise non-isomorphic simple $R$-modules $E_i$ and $k,n_i \geq 1$.
Set $P_i:=\Ann{R}{E_i}$. Then $S=\oplus_{i=1}^k
M_{n_i}(\Delta_i)$ where $\Delta_i=\EndX{R}{E_i}$. Assume that $R$
is a PI-ring. By Kaplansky's Theorem (see \cite[13.3.8]{McConnellRobson}) there exists $m_i\geq 1$ such that
$R/P_i$ is isomorphic to the full matrix ring $M_{m_i}(\Delta_i)$ and
$E_i$ is a finite-dimensional $\Delta_i$-vector space. Hence $E_i^{n_i}$ and also $M$ are finitely generated over their
endomorphism rings.
\end{proof}

\subsection{} Applying the above Lemma to the module algebra
situation we will see, that the $H$-action on a semiprime Goldie
PI module algebra whose non-zero $H$-stable ideals contain
non-zero central $H$-invariant elements can be extended to its
ring of quotients.

\begin{prop}\label{ExistenceHInvariants} Let $H$ be a Hopf algebra over $R$ with $H_R$ finitely
generated and let $A$ be a semiprime Goldie PI $H$-module algebra
with classical ring of quotients $\Ql{A}$. If every non-zero
$H$-stable ideal of $A$ contains a non-zero central $H$-invariant
element, then the $H$-action on $A$ can be extended to $\Ql{A}$
and $\Ql{A}$ is equal to the central localisation $A[\cC^{-1}]$ of
regular elements $\cC$ of the subring $Z(A)^H$ of $A$.
\end{prop}

\begin{proof}
Let $Z(A)^H:=Z(A)\cap A^H$ and let $\cC$ denote the set of regular
elements of $Z(A)^H$. The elements of $\cC$ form an Ore set in $A$
and are also regular elements of $A$ since $\Ann{A}{c}^H =0$
implies $\Ann{A}{c}=0$ for all $c\in \cC$. Denote by
$\tA:=A[\cC^{-1}]$ the localisation of $A$ by $\cC$. Note that $A$
is a subring of $\tA$ and the map $I\mapsto (I\cap A)$ from ideals
of $\tA$ to ideals of $A$ is injective. In particular $\tA$ is
semiprime. Since $\tA$ is a central extension of the PI-ring $A$,
$\tA$ is PI by \cite[13.1.11]{McConnellRobson}. By
\cite[6.1.1]{Rowen} $\tA\otimes \tA^{op}$ is a PI-ring and hence
its factor ring $$\tA\otimes \tA^{op}/\Ann{\tA\otimes \tA^{op}}{\tA} \simeq M(\tA) := \langle \{L_x, R_x \mid x\in\tA
\}\rangle \subseteq \EndX{R}{\tA}$$ is a PI-ring. The $H$-action
on $A$ extends trivially to $\tA$ by letting an element $h\in H$
act on an element $ac^{-1}$ as $(h\cdot a)c^{-1}$. Since $H$ is
finitely generated as $R$-module, $M_H(\tA)$ is a finite extension
of $M(\tA)$ and therefore also a PI-ring by
\cite[13.4.9]{McConnellRobson}. Note that
$$\EndX{M_H(\tA)}{\tA} \simeq Z(\tA)^H \simeq Z(A)^H[\cC^{-1}]\simeq \Ql{Z(A)^H}$$
is semisimple artinian.
Moreover let $I$ be a non-trivial $H$-stable ideal of $\tA$. Then
$I\cap A$ is a non-trivial $H$-stable ideal of $A$ and contains a
non-trivial central $H$-invariant element. Using the isomorphism
$$\HomX{M_H(\tA)}{\tA}{I} \simeq I\cap Z(\tA)^H\neq 0$$
we see that $\tA$ is a retractable module over the PI-ring $M_H(\tA)$ having a semisimple artinian endomorphism ring
isomorphic to $Z(\tA)^H$.
By Lemma \ref{RetractableSemisimple}, $\tA$ is finitely generated
over $Z(\tA)^H$ and is therefore left and
right artinian. Being semiprime artinian makes $\tA$ a semisimple
artinian ring and since $A$ is a left order in $\tA$ we can
conclude that $\tA$ is equal to the classical ring of quotients of
$A$. Thus $\Ql{A}=\tA$ is finitely generated over $Z(\Ql{A})^H$.
\end{proof}

\subsection{}
In case there do not exist non-trivial $H$-stable ideals we obtain the
following corollary from the previous proposition.

\begin{cor}\label{CorExistenceHInvariants}
Let $H$ be a Hopf algebra over $R$ with $H_R$ finitely generated.
Any semiprime Goldie PI $H$-module algebra that is $H$-simple is finite
dimensional over $Z(A)^H$ and equals its classical ring of quotients.
\end{cor}

\begin{proof} Since $A$ is $H$-simple $Z(A)^H$ is a field. Thus by Proposition \ref{ExistenceHInvariants}
$\Ql{A}=A$ and $dim_{Z(A)^H}(A)$ is finite.
\end{proof}

\subsection{} We can now prove the main result of this section
showing that the ability of extending the $H$-action to the classical left ring
of quotients of a semiprime Goldie PI $H$-module algebra $A$ with
central invariants is equivalent to $\AH$ being semiprime.

\begin{thm}\label{SemiprimeGoldieZentral} Let $H$ be a Hopf algebra over $R$ with
$H_R$ finitely generated and projective.
Let $A$ be a semiprime Goldie PI $H$-module algebra with central
invariants such that there exists a left or right integral $t$
with $\eps(t)1_A$ invertible in $A$. Then the following statements
are equivalent:
\begin{enumerate}
\item[(a)] Every essential left ideal of $A$ contains a regular $H$-invariant element.
\item[(b)]  The $H$-action on $A$ extends to the classical left ring of quotients $Q_{cl}(A)$.
\item[(c)] $\AH$ is semiprime.
\item[(d)] Every $H$-stable left ideal of $A$ contains a non-zero $H$-invariant element.
\end{enumerate}
Then $\Ql{A}=A[\cC^{-1}]$ and $\Ql{\AH}=\AH[\cC^{-1}\# 1]$, where $\cC$ denotes the set of regular elements of $A^H$.
\end{thm}

\begin{proof}
Let $\cC$ denote the set of regular elements of $A^H$.\\
$(a)\Rightarrow (b)$ Consider $\tA:=A[\cC^{-1}]$ and let
$I$ be an essential left ideal of $\tA$.
Then $I\cap A$ is an essential left ideal of $A$ and contains an element of $\cC$. Hence $I=\tA$ shows that
$\tA$ has no proper essential submodules and must be semisimple artinian. Since $A$ is a right order in $\tA$ we
obtain that $\tA=\Ql{A}$. The $H$-action can be extended trivially to $\tA$.\\
$(b)\Rightarrow (c)$ Let $\cD$ denote the set of regular elements
of $A$. The $H$-action on $A$ can be extended to the classical
left ring of quotients $\Ql{A}=A[\cD^{-1}]$ by hypothesis. Since $A$
is a semiprime Goldie PI-algebra, $\Ql{A}$ is semisimple artinian.
By Corollary \ref{VonNeumannRegularesSmashProdukt} $\Ql{A}\# H$ is semisimple
artinian since $\eps(t)1_{\Ql{A}}$ is invertible in $\Ql{A}$.
As $A$ is a left and right order in $\Ql{A}$ every element
of $\Ql{A}$ can be written in the form $d^{-1}a$ with $d\in \cD$
and $a\in A$. Hence $\AH$ is a left order in $\Ql{A}\# H$.  Thus
by Goldie's Theorem $\AH$ is semiprime and $\Ql{\AH} \simeq
\Ql{A}\# H$.\\
$(c)\Rightarrow (d)$ Note that $a \mapsto a\# t$ is an injective
$\AH$-linear map from $A$ to $\AH$. Assume that  $\AH$ is
semiprime and let $I$ be a non-zero $H$-stable left ideal of $A$.
Then $0\neq (I\#t)^2 = I(t\cdot I)\# t$ shows $I^H \supseteq t\cdot I \neq 0$.\\
$(d)\Rightarrow (a)$
By Proposition \ref{ExistenceHInvariants} $\tA=A[\cC^{-1}]$ equals $\Ql{A}$ and
is semisimple artinian.
Let $I$ be an essential left ideal of $A$. Then $I[\cC^{-1}]$ is an essential left
ideal of the semisimple ring $\tA$ and therefore improper. Thus $I[\cC^{-1}]=\tA$ implies that there exist
$a\in I$ and $c\in\cC$ such that $ac^{-1}=1$. Equivalently $a=c\in I\cap \cC$ shows
that $I$ contains a regular $H$-invariant element.
\end{proof}

\subsection{} Note that condition $(d)$ of
\ref{SemiprimeGoldieZentral} says that for every left ideal $I$ in
the filter $\cF$ of essential left ideals of $A$ and for every
$h\in H$ there exists an essential left ideal $I' \in \cF$ such
that $hI' \subseteq I$. Montgomery had termed $H$-actions with
this property $\cF$-continuous and had shown in \cite{Montgomery}
that this condition is sufficient for extending the $H$-action to
the ring of quotients with respect to the filter $\cF$. We see
that under the assumptions of \ref{SemiprimeGoldieZentral} the
$\cF$-continuity of the $H$-action is also a necessary condition.
\subsection{}
Combining Theorem \ref{CommutativeCase} and Theorem \ref{SemiprimeGoldieZentral} we obtain the following Corollary for
Hopf actions on integral domains.

\begin{cor} Let $H$ be a semisimple, cosemisimple Hopf algebra over a field $k$ and let $A$ be
a left $H$-module algebra that is an integral domain. Then the quotient field $Q$ of $A$ equals $A[\cC^{-1}]$ where
$\cC:= A^H \setminus \{0\}$. The $H$-action extends to $Q$ and $Q^H\subseteq Q$ is a finite field extension.
$\AH$ is a semiprime Goldie PI-algebra with classical ring of quotient isomorphic to $Q\# H$.
\end{cor}

\subsection{}
A classical result of Bergman and Isaac asserts, that a ring $A$
with group action $G$ such $A$ is $|G|$-torsionfree is nilpotent
whenever $A^G$ is nilpotent. As a kind of Hopf-algebraic analogue
Bahturin and Linchenko showed in \cite{BahturinLinchenko} that
every left $H$-module algebra $A$ (possibly without unit) is
nilpotent whenever $A^H$ is nilpotent if and only if every left
$H$-module algebra $A$ (possibly without unit) is PI whenever
$A^H$ is PI if and only if $T(H)/\langle\int_l\rangle$ has finite
dimension, where $H$ is a finite dimensional Hopf algebra over a
field of characteristic $0$, $T(H)$ denotes the tensor algebra of
$H$ and $\langle\int_l\rangle$ the
 ideal of $T(H)$ generated by the left integrals in $H$.
They also show that under those equivalent conditions above $H$
must be semisimple. Whether every semisimple Hopf algebra fulfills
one of the above properties is still open.

Combining Bahturin and Linchenko's result with Theorem
\ref{SemiprimeGoldieZentral} we can conclude the following: If $H$
is a finite dimensional Hopf algebra over a field $k$ of
characteristic $0$ such that $T(H)/\langle \int_l\rangle$ is
finite dimensional and if $A$ is a  semiprime Goldie left
$H$-module algebra with central invariants then one can extended
the $H$-action to $\Ql{A}$, $\Ql{A}$ is equal to the localisation
of $A$ by the regular elements of $A^H$ and $\AH$ is semiprime
with classical ring of quotients equal to $\Ql{A}\# H$.

%\begin{proof} By Bahturin and Linchenkos theorem, $A$ is an
%PI-algebra. Moreover if $I$ is any $H$-stable left ideal with
%$I^H=0$, then by Bahturin and Linchenko's theorem $I$ is
%nilpotent. But $A$ being semiprime implies $I=0$. Hence property
%$(d)$ of Theorem \ref{SemiprimeGoldieZentral} is fulfilled and the
%result follows.
%\end{proof}

\section{Drinfeld Twists of strongly semisimple Hopf algebras}

We finish the paper by showing that Cohen's question has a
positive answer if $H$ is semisimple cosemisimple triangular.

\begin{defn}
A Hopf algebra $H$ over $R$ is called  {\it strongly semisimple } if
for every $H$-semiprime left $H$-module algebra $A$ the smash
product $\AH$ is semiprime.
\end{defn}

Criterions for a Hopf algebra to be strongly semisimple are given
in  \cite{MontgomerySchneider99} but those criterions are hard to verify.
Over a field, every
commutative or cocommutative semisimple Hopf algebra is strongly
semisimple. Moreover Montgomery and Schneider showed, that every
semisimple Hopf algebra that admits a
normal series $H_i$ whose quotients $H_{i+1}/H_i$ are either
commutative or cocommutative, is strongly semisimple (see
\cite[8.16]{MontgomerySchneider99}). Those Hopf algebras are called
semi-solvable.

We will show that the class of strongly semisimple Hopf algebras
is closed under Drinfeld twists. Applying a theorem of Etingof and
Gelaki, that classifies all triangular semisimple cosemisimple
Hopf algebras as Drinfeld Twists of group algebras, we can
conclude that all triangular semisimple cosemisimple Hopf algebras
are strongly semisimple.

\subsection{} Recall the definition of Drinfeld twists for a Hopf algebra.
\begin{defn}
Let $H$ be an Hopf algebra over $R$. A {\it Drinfeld
Twist}\index{Drinfeld twist} for $H$ is an invertible element
$J\in H\otimes H$, such that
\begin{equation}\label{Twist1}
 (J\otimes 1)(\Delta \otimes 1)(J)  = (1\otimes J)(1\otimes \Delta)(J)
\end{equation}
\begin{equation}\label{Twist2}
(\eps\otimes 1)(J) = 1 = (1 \otimes \eps)(J)
\end{equation}
holds. We write formally $J=:\sum J^1 \otimes J^2$ and
$J^{-1}=:Q=:\sum Q^1 \otimes Q^2$.
\end{defn}

If $H$ is a Hopf algebra over $R$ with comultiplication $\Delta$ and antipode $S$, then
$\Delta^J := J\Delta J^{-1}$ defines a new comultiplication on $H$ with
$\Delta^J(h) := J\Delta(h)J^{-1}$ for all $h \in H$.
Let $U:=\sum J^1S(J^2)$ and $U^{-1}=\sum S(Q^1)Q^2$ and define a new map
$S^J := USU^{-1}$ by $S^J(h):= US(h)U^{-1}$ for all $h\in H$. Then it has been shown in
\cite[2.3.4]{MajidBook} that $\Delta^J$ and $S^J$ define a new Hopf algebra structure on $H$ keeping the same multiplication, unit and counit.
We denote the obtained Hopf algebra by $H^J$. Obviously $\Delta^J(h) J = J \Delta(h)$ for all $h\in H$.

%\subsection{}\label{DrinfeldTwistLemma}
Moreover it is not difficult to see that $J^{-1}$ is a Drinfeld twist for $H^J$.
%\begin{lem}
%Let $J$ be a Drinfeld twist for a Hopf algebra $H$ over $R$. Then
%$J^{-1}$ is a Drinfeld twist for $H^J$.
%\end{lem}

%\begin{proof}
%Note that $(\Delta \otimes 1)(J^{-1})$ is an inverse for $(\Delta
%\otimes 1)(J)$ since $\Delta$ is a homomorphism of algebras. We
%reformulate equation (\ref{Twist1}):
%\begin{eqnarray*}
%                & (J\otimes 1)(\Delta \otimes 1)(J)  = (1\otimes J)(1\otimes \Delta)(J)\\
%\Leftrightarrow & 1 \otimes 1 \otimes 1  = (1\otimes J)(1\otimes \Delta)(J)(\Delta \otimes 1)(J^{-1})(J^{-1}\otimes 1)\\
%\Leftrightarrow & (1\otimes \Delta)(J^{-1})(1\otimes J^{-1})  = (\Delta \otimes 1)(J^{-1})(J^{-1}\otimes 1)\\
%\Leftrightarrow & (1\otimes J^{-1})(1\otimes \Delta^J)(J^{-1})  =
%(J^{-1}\otimes 1)(\Delta^J \otimes 1)(J^{-1})
%\end{eqnarray*}
%where last step follows from $J^{-1}\Delta^J = \Delta J^{-1}$.
%Hence equation (\ref{Twist1}) holds for $J^{-1}$ in $H^J$.
%Equation (\ref{Twist2}) holds, since $\eps$ is a homomorphism of
%algebras.
%\end{proof}

\subsection{}\label{DefinitionTwistedModuleAlgebra}
Having `twisted' the comultiplication of $H$ we can also
`twist' the multiplication of a left $H$-module algebra $A$ such that $A$ becomes a left $H^J$-module algebra.

\begin{defn}
Let $A$ be a left $H$-module algebra with multiplication $\mu$ and
let $J$ be a Drinfeld twist for $H$. We define a new
multiplication $\mu^J: A\otimes A \lra A$ on $A$
with
$$a \cdot_J b := \mu^J(a\otimes b) := \sum (Q^1\cdot a)(Q^2\cdot b) \:\:\mbox{ for all }a,b \in A.$$
\end{defn}

It had been shown in \cite[2.3.8]{MajidBook} that $A^J$ with multiplication $\mu^J$ is a left $H^J$-module algebra.
Moreover the smash products $\AH$ and $A^J \# H^J$ are isomorphic $R$-algebras.
This follows from a more general theorem by Majid (see  \cite[2.9]{Majid}).

\subsection{}\label{DrinfeldTwistProp} Note that for every two elements $a,b \in A$ we have:
$$ a b =\sum (J^1\cdot a)\cdot_J (J^2\cdot b).$$
In particular take any $H^J$-stable ideal $I$ of $A^J$, then $I$ is also an $H$-stable ideal of $A$.
Moreover if $I$ is nilpotent as an ideal of $A^J$, then it is also nilpotent as an ideal of $A$.
This shows $A^J$ is $H^J$-semiprime whenever $A$ is $H$-semiprime.
By the same argument applied to $A=(A^J)^{J^{-1}}$ we obtain $A$ is $H$-semiprime whenever $A^J$
is $H^J$-semiprime.

\subsection{}\label{AbschlussTwists} Combining the results of the last two paragraphs
%\ref{DrinfeldTwistLemma} and \ref{DrinfeldTwistProp}
we can prove that the class of strongly semisimple Hopf algebras is closed under Drinfeld twists.
\begin{cor}
The class of strongly semisimple Hopf algebras is closed under
Drinfeld twists.
\end{cor}
\begin{proof}
Let $H$ be a strongly semisimple Hopf algebra and let $J$ be a Drinfeld twist for $H$.
Let $A$ be a left $H^J$-module algebra, then $A^{J^{-1}}$ is a
left $H$-module algebra by \cite[2.3.8]{MajidBook}. If $A$ is
$H^J$-semiprime, then $A^{J^{-1}}$ is $H$-semiprime by
\ref{DrinfeldTwistProp}. As noticed in \ref{DefinitionTwistedModuleAlgebra} from
 \cite[2.9]{Majid} follows $$A^{J^{-1}} \# H = A^{J^{-1}} \# {H^J}^{J^{-1}} \simeq A\# H^J.$$
Since $H$ is strongly semisimple, $A^{J^{-1}}\# H$ and therefore
$A\# H^J$ is semiprime. Hence $H^J$ is strongly semisimple.
\end{proof}

%\subsection{} A Hopf algbera is called {\it simple} if it does not contain any
%non-trivial normal sub-Hopf algebra. The only known simple
%semisimple Hopf algebras over a field $k$ are group rings $k[G]$
%of finite simple groups $G$ with $char(k)\nmid |G|$ their
%Drinfeld twists $k[G]^J$ and their duals $k[G]^\ast$ .
%Since $k[G]$ is strongly semisimple any Drinfeld twist $k[G]$ is also strongly semisimple by
%Corollary \ref{AbschlussTwists}. The fact that $k[G]^\ast$ is strongly semisimple is due to Cohen and Rowen.
%Hence we proved that all `known' simple semisimple Hopf algebras are strongly semisimple.
%Whether there are other simple semisimple Hopf algebras is an open question (see \cite[Question
%2.3]{Andruskiewitsch}).

\subsection{}A Hopf algebra is called {\it triangular}, if there exists an
invertible element $\cR \in H\otimes H$ with
$$ (\Delta \otimes 1)(\cR) = \cR_{13}\cR_{23},
\: (1\otimes \Delta)(\cR)=\cR_{13}\cR_{12}, \: \Delta^{cop} =
\cR\Delta\cR^{-1} \mbox{ und } \cR^{-1}=\tau(\cR)$$ where
$\tau:H\otimes H\lra H\otimes H$ is the isomorphism $x\otimes
y\longmapsto y\otimes x$. For $\cR= \sum a_i \otimes b_i$ we set
$$\cR_{13}:= \sum a_i \otimes 1 \otimes b_i,
\cR_{23}:= \sum 1\otimes a_i  \otimes b_i, \cR_{12}:= \sum a_i
\otimes b_i \otimes 1.$$

 P.Etingof and S.Gelaki classified in \cite{EtingofGelaki00}
 semisimple cosemisimple triangular Hopf algebras
 over algebraically closed fields as Drinfeld twists of group rings.
 From this we obtain as a corollary:

\begin{cor}
All triangular semisimple cosemisimple Hopf algebras over an
algebraically closed field are strongly semisimple.
\end{cor}

\begin{proof}
Let $H$ be a semisimple cosemisimple triangular Hopf algebra over
an algebraically closed field $k$. By Etingof and Gelaki's result
\cite[Corollary 6.2]{EtingofGelaki00} there exists a group $G$ and
a Drinfeld twist $J \in k[G] \otimes k[G]$ such that $H\simeq
k[G]^J$ as Hopf algebras. As $k[G]$ is strongly semisimple also
$H$ is strongly semisimple by \ref{AbschlussTwists}.
\end{proof}

\section*{Acknowledgment}

Parts of the material presented here are included in the author's doctoral thesis at
the Heinrich-Heine Universit\"at D\"usseldorf. The author would like to express his gratitude to his supervisor
Professor Robert Wisbauer for all his help advice and encouragement. Moreover he would like to thank the
Mathematics Department of the University of Porto for suspending him from teaching duty while writing his thesis.
The author would also like to thank Professor Kostia Beidar for helpful discussions on the topic and
Professor Shahn Majid for pointing out reference \cite[2.9]{Majid}. Last but not least the author would like to express
his gratitude to his wife Paula for all her patience and love during this time.

\bibliographystyle{amsplain}

%\bibliography{promo}
\end{document}